\documentclass[11pt]{article}

\usepackage[english,greek]{babel}
\usepackage[utf8x]{inputenc}

\usepackage{amsmath,amsthm,amssymb,stmaryrd,bigints,mathabx}

\usepackage[pdftex]{graphicx}
\usepackage{subcaption}
\usepackage[titletoc,title]{appendix}
\usepackage{datetime,color,cancel,listings,authblk}
\usepackage[hyperfootnotes=false]{hyperref}
\hypersetup{
  colorlinks = true,
  linkcolor = black,
  citecolor = blue,
  urlcolor = black
}
\usepackage[margin=2.0cm]{geometry}

\usepackage{float}

\floatstyle{ruled}
\newfloat{algorithm}{thp}{lop}
\floatname{algorithm}{Algortithm}

\newcommand{\E}{\mathbb{E}}

\newcommand{\R}{\mathbb{R}}

\newcommand{\C}{\mathbb{C}}
\newcommand{\beq}{\begin{equation}}
\newcommand{\ee}{\end{equation}}
\newcommand{\bac}{\begin{array}{c}}
\newcommand{\ea}{\end{array}}
\newcommand{\bal}{\begin{aligned}}
\newcommand{\eal}{\end{aligned}}

\newcommand{\real}{\operatorname{Re}}
\newcommand{\imag}{\operatorname{Im}}

\newcommand{\BFI}{\textrm{BFI}}

\newcommand{\Expe}{\operatorname{E}}

\newcommand{\vertiii}[1]{{\left\vert\kern-0.25ex\left\vert\kern-0.25ex\left\vert #1 
    \right\vert\kern-0.25ex\right\vert\kern-0.25ex\right\vert}}

\begin{document}
\newtheorem{theorem}{Theorem}[section]
\newtheorem{lemma}[theorem]{Lemma}
\newtheorem{remark}[theorem]{Remark}
\newtheorem{observation}[theorem]{Observation}
\newtheorem{definition}[theorem]{Definition}
\newtheorem{example}[theorem]{Example}
\newtheorem{corollary}[theorem]{Corollary}
\newtheorem{assumption}{Assumption}
\newtheorem{property}{Property}

\selectlanguage{english}

\title{Modelling of ocean waves with the Alber equation: \\application to non-parametric spectra and generalization to crossing seas}

\author[1]{Agissilaos G. Athanassoulis}
\author[2]{Odin Gramstad}

\affil[1]{Department of Mathematics, University of Dundee}
\affil[2]{Hydrodynamics, MetOcean \& SRA, Energy Systems, DNV}

\maketitle

\begin{abstract} The Alber equation is a phase-averaged second-moment model for the statistics of a sea state, which recently has been attracting renewed attention. We extend it in two ways: firstly, we derive a generalized Alber system starting from a system of nonlinear Schr\"odinger equations, which contains the classical Alber equation as a special case but can also describe crossing seas, i.e. two wavesystems with different wavenumbers crossing. (These can be two completely independent wavenumbers, i.e. in general different directions and different moduli.) We also derive the  associated 2-dimensional scalar instability condition. This is the first time that a modulation instability condition applicable to crossing seas  has been systematically derived for general spectra. Secondly,  we use the classical Alber equation and its associated instability condition to quantify how close a given non-parametric spectrum is to being modulationally unstable. We apply this to a dataset of 100 non-parametric spectra provided by the Norwegian Meteorological Institute, and find the vast majority of realistic spectra turn out to be stable, but three  extreme sea states are found to be unstable (out of 20 sea states chosen for their severity). Moreover, we introduce a  novel ``proximity to instability'' (PTI) metric, inspired by the stability analysis. This is seen to correlate strongly with the steepness and Benjamin-Feir Index (BFI) for the sea states in our dataset ($>85\%$ Spearman rank correlation). Furthermore, upon comparing with phase-resolved broadband Monte Carlo simulations, the kurtosis and probability of rogue waves for each sea state are also seen to correlate well with its PTI ($>85\%$ Spearman rank correlation).

\end{abstract}

\noindent {\bf Keywords: } ocean waves, Alber equation, modulation instability, rogue waves, crossing seas

\setcounter{tocdepth}{1}
\tableofcontents

\section{Introduction}

Ocean waves are a very active field of mathematical modelling and analysis. The first principles hydrodynamic models for gravity waves are by now well understood
 \cite{Alazard2014,Lannes2005,Luke1967}. An array of approximate models is well established and widely used, including the nonlinear Schr\"odinger equation (NLS) and its variants \cite{benney1967propagation,Trulsen1996}, the Zakharov equation \cite{zakharov1968stability}, the coupled-mode systems  \cite{Athanassoulis1999,Athanassoulis2017a,Papoutsellis2018}, the High Order Spectral Method (HOSM) \cite{dommermuth1987high,west1987new,gouin2016development} and others. (In shallow water an even larger collection of models is being used, but here we focus on deep water.) One reason for the wide use of approximate models in oceanography is the need to study large wavefields, with hundreds or thousands of individual wavelengths. Thus there is a trade-off between hydrodynamic fidelity and the ability to scale up models, as well as the accessibility of powerful qualitative insights.
 
 In fact, actual ocean waves include multi-physics and non-differentiable  phenomena, such as wind forcing, wave breaking etc. These are not accounted for in the classical ``exact'' hydrodynamic model anyway, and are still being understood on different levels \cite{Brunetti2014a,Brunetti2014,Armaroli2018,Buhler2016,Janssen}.

In this paper we will use {stochastic modelling} of ocean waves and, furthermore, explore {what phase-averaged stochastic models may reveal about rogue waves in particular}. This  will lead us to novel mathematical results;  it is also worth mentioning that  this kind of overall approach has been identified as a priority within the broader marine research \& industry community \cite{15PrioritiesforWindWavesResearchAnAustralianPerspective}.

\medskip


The most well-known stochastic models for ocean waves include the CSY equation \cite{Crawford1980,Andrade2020a,Andrade2020b}, Hasselmann's equation \cite{hasselmann1962non} and the Alber equation \cite{Alber1978} that we will focus on here. For a recent review of various stochastic models one can see \cite{stuhlmeier2019nonlinear}. Broadly speaking, they are {moment equations}, starting from {phase-resolved equations} for the sea surface (such as Zakharov's equation or the  NLS) as an approximation for deterministic wave dynamics. One then takes stochastic moments of the deterministic equations; due to the nonlinearity of these equations, an infinite hierarchy of moments is produced. A {gaussian second-order moment closure} is then used, to produce a closed equation for the second stochastic moment. The resulting equation is  {phase-averaged}, meaning that in no longer resolves individual wave peaks and troughs, but instead the evolution and propagation of the statistics of the wavefield.

A key approximation step is the gaussian moment closure. This is of course not exact, but in many cases the free surface is indeed close to being gaussian \cite{Komen1994,Ochi1998}, making the gaussian closure plausible. It is important to keep in mind that fidelity to the deterministic, phase-resolved model is not the only consideration: for example, the  exact  infinite hierarchy of moments, is not just more complicated as a mathematical model: {it is impossible to initialize meaningfully}. The vast majority of synoptic data collected from the ocean is, or can be converted to, some kind of second moment \cite{Ochi1998}. There is some data involving moderately higher moments, but very little data or know-how exists for moments higher than 4$^{th}$ order.

The question of using realistic data for initialisation is extremely important; all the more so  in the study of extreme sea states and rogue waves. For example, {every deterministic model can produce  ``rogue waves'' on demand, by carefully preparing particular  initial conditions;} however, the real-life question is {how often would these ``initial conditions leading to rogue waves'' realistically appear?} A stochastic approach can directly attack this question e.g. by a phase-resolved Monte Carlo  approach \cite{Toffoli2010}; this has a number of advantages, but is clearly expensive to be applied indiscriminantly. Another possibility would be using phase averaged stochastic models, like the moment equations discussed above, to directly investigate whether a sea state is likely to support the rapid concentrations of energy. That way sea states of interest could be selected and computational resources focused on them. In fact, it appears that the sea states highlighted as more unstable by the Alber equation do turn out to exhibit a higher probability of extreme events in a phase-resolved  Monte Carlo simulation (details in Section \ref{sec:summary} and Figure \ref{fig:82}). 

%
%

\section{Ocean wave modelling with the Alber equation}

The Alber equation, its derivation and interpretation have been widely studied and explained. We refer the interested reader to \cite{Alber1978,Ribal2013a,Gramstad2017,Athanassoulis2018} and the references therein for more details. Here we will briefly present its derivation and main features in order to make the paper self-contained.
 
\subsection{Derivation}

The 
cubic focusing nonlinear Schr\"odinger equation (NLS)
\beq\label{eq:nls}
i \partial_t u - p \Delta u - q |u|^2 u =0.
\ee
 is an approximate model for the {envelope} of a narrow-banded wavetrain with {carrier wavenumber $\boldsymbol{k_0}$} along its direction of propagation (unidirectional propagation), cf. e.g. \cite{Chiang2005}.
Thus, in deep water the sea surface elevation $\eta(x,t)$ is related to the complex-valued envelope $u(x,t)$ through
\[
\eta(x,t) = \real \left[ u(x,t) e^{i(k_0 x - \omega_0 t)}\right], \qquad \omega_0 = \sqrt{g k_0}.
\]
The coefficients of the equation also depend on $k_0=|\boldsymbol{k_0}|,$
\[
p=\frac{\sqrt{g}}{8 k_0^{\frac{3}2}}, \qquad q=\frac{\sqrt{g}}2 k_0^{\frac{5}2}.
\]
This is an  asymptotic model, where the order parameter is the steepness of the waves, and moreover assumes narrow-bandedness of the wavefield around the carrier wavenumber. The NLS and its variants are widely used as relevant to a wide array of realistic scenarios \cite{Toffoli2010}.\footnote{There are some key facts here; water waves cannot get too steep, since they break at a slope of about $0.13,$ and typically they are quite less steep than that. Moreover,  wavelengths for gravity waves vary over no more than two orders of magnitudes ($\sim 5-500m$). In fact, sea-state spectra are known to often   be more narrowly supported  even within that range. 
In contrast, {crossing seas} (i.e. $O(1)$ wave-systems coming from different directions) cannot be considered narrowband due to the different directions; thus, they provide a prime example of a realistic situation where the NLS \eqref{eq:nls} would {not} be a satisfactory model. }
A different kind of limitation of the NLS is that it is a {deterministic phase-resolved model; i.e. any prediction with it will not be better or more accurate than the initial condition used.} But realistic wave-systems are not widely available as phase-resolved initial conditions, as opposed to power spectra. In that context, Alber \cite{Alber1978} proposed generating a second-order moment equation from the NLS, considered with stochastic initial data. 
Denoting by 
\[
R(x,y,t)= \E[u(x,t)\overline{u}(y,t)]
\]
the autocorrelation of the envelope $u,$
one obtains
\beq\label{eq:preclo}
i \partial_{t} R(x,y,t)+{p}\left(\Delta_{x}-\Delta_{y}\right) R(x,y,t)+q \E \Big[u(x,t)\overline{u}(y,t)[u(x,t)\overline{u}(x,t)-u(y,t)\overline{u}(y,t)]\Big]=0.
\ee
By using the gaussian closure
\beq\label{eq:gmcl098}
\Expe[|u(\alpha,t)|^{2} u(\alpha,t)\overline{u}(\beta,t)] = 2 R(\alpha,\alpha,t)R(\alpha,\beta,t),
\ee
the autocorrelation can now be seen to satisfy
\beq\label{eq:Alber1111}
i \partial_{t} R(x,y,t)+{p}\left(\Delta_{x}-\Delta_{y}\right) R(x,y,t)+2q R(x, y, t)[R(x, x, t)-R(y, y,t)]=0.
\ee

A key restriction the Alber equation inherits from its starting point, the  NLS equation,  is {narrowbandedness}. This is arguably not too restrictive for unidirectional sea states \cite{Toffoli2010}; however it completely fails for crossing seas, i.e. sea states where several different directions carry substantial wave energy. To address this limitation, we will derive a generalized Alber equation valid in a crossing seas scenario. This is made possible by starting from a system of NLS equations describing the crossing wavetrains which was derived in  \cite{Hammack2005}.  The resulting generalized Alber equation is reported in Section \ref{sec:mainresults}, and derived in detail in Section \ref{sec:crossseas}.  Crossing seas are receiving increased attention as a possible incubator of rogue waves \cite{Onorato2010,Gramstad2018}, so it is significant to have a moment equation applicable to such scenarios.

\subsection{The stability-of-homogeneity question} 

It is empirically known that sea-states are typically homogeneous and stationary, at least for appropriate lengthscales and  timescales \cite{Komen1994,Ochi1998,Labeyrie1990,Tournadre1993}. This feature is reflected in the Alber equation; for example observe that $R(x,y,t)=\Gamma(x-y)$ is a solution of equation \eqref{eq:Alber1111} for any smooth function $\Gamma.$ The next step is to investigate the stability of such homogeneous and stationary solutions:
 let us consider a {weakly  inhomogeneous initial sea-state}, i.e. assume that the autocorrelation is initially of the form
\beq\label{eq:weaklynonhomogIC}
R(x,y,0)= \E[u(x,0)\overline{u}(y,0)] = \Gamma(x-y) + \epsilon \rho(x,y,0), \qquad \epsilon \ll 1
\ee
for some nice functions $\Gamma,\rho$ (e.g. Schwartz class test functions).   By inserting eq. \eqref{eq:weaklynonhomogIC} in eq. \eqref{eq:Alber1111} one can reformulate the problem in terms of the inhomogeneity $\rho,$
\beq\label{eq:Alber11112}
i \partial_{t} \rho(x,y,t)+{p}\left(\Delta_{x}-\Delta_{y}\right) \rho(x,y,t)+2q [\Gamma(x-y) + \epsilon\rho(x,y,t)][\rho(x, x, t)-\rho(y, y,t)]=0.
\ee


So stability for homogeneous state-states is controlled by the boundedness (or lack thereof) of the inhomogeneity $\rho$ in equation \eqref{eq:Alber11112}. Can $\rho$ grow in time to the extent that $\epsilon \rho(x,y,t)$ is no longer small? Or is there a guarantee that $\rho$ stays bounded? In the latter case of stability, the autocorrelation will simply stay close to $\Gamma(x-y)$ for all times. In the unstable case however, even if initially very close to homogeneous, the autocorrelation could develop significant inhomogeneities. 

\begin{figure}[!h]
\begin{center}
\includegraphics[width=0.9\textwidth]{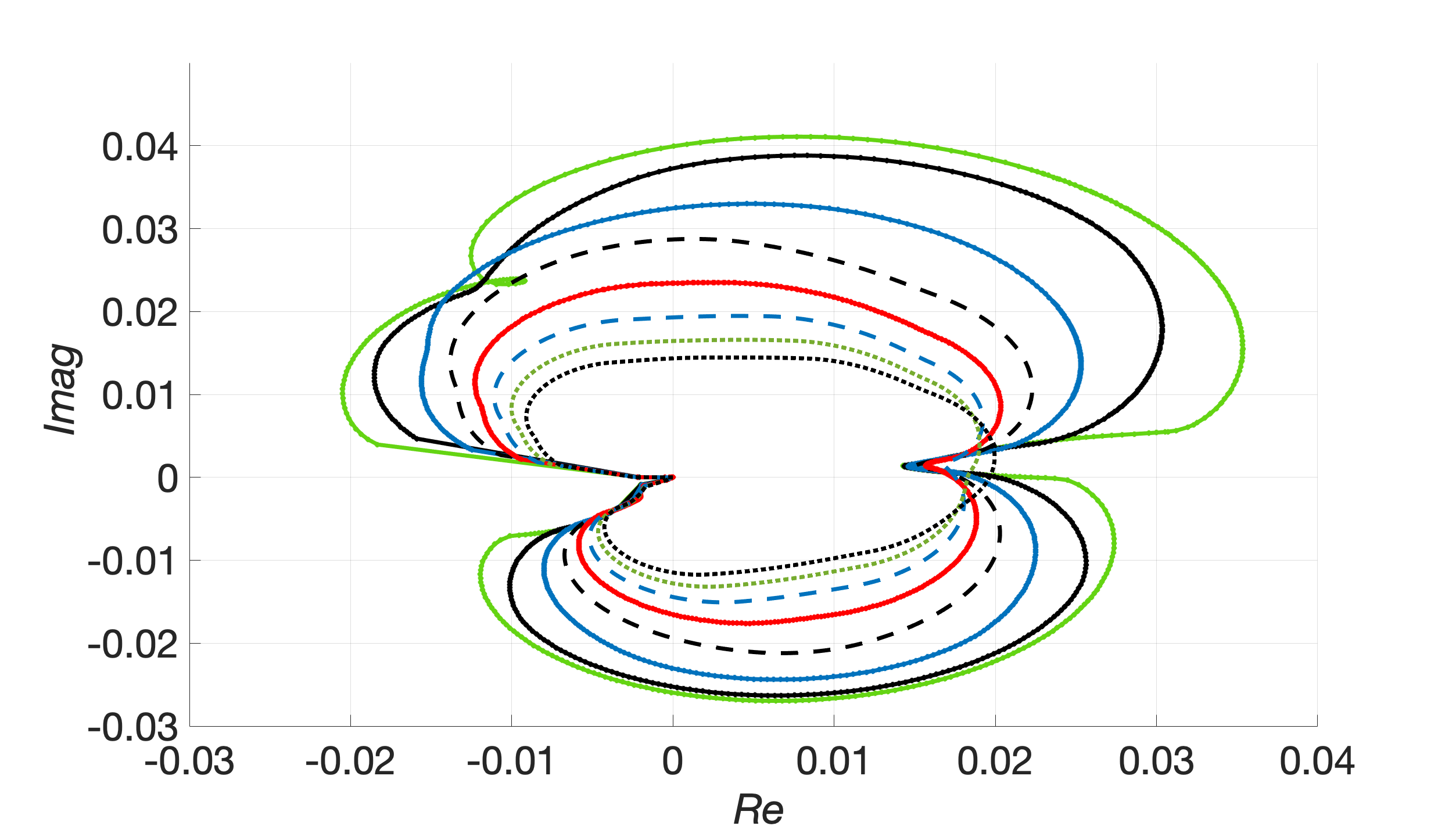} 
\caption{
Curves $\Gamma_X$ for a nonparametric spectrum  on the complex plane.  The smaller closed curves correspond to larger $|X|;$ here $X=(1+120\cdot n)5\cdot 10^{-4},$ $n=0,\dots,7.$ Taking $|X|<5\cdot 10^{-4}$ does not change the outermost curve noticeably, as $D_XP$ has effectively converged to $P'.$ 
The real number $1/4\pi \approx 0.08$ is {not} circumscribed by any of the curves, i.e. the spectrum does {not} exhibit modulation instability (cf. Definition \ref{def:penralbcond}).  The spectrum used is a unimodal spectrum with a reference steepness of around $7\%,$ taken out of the dataset of 100 nonparametric spectra.
}
\label{fig:1}
\end{center}
\end{figure}

In \cite{Alber1978} a sufficient condition for linear instability was derived in terms of the {spectrum} $S(k),$ namely the Fourier transform of the autocorrelation function $\Gamma,$ 
\[
S(k)=\mathcal{F}_{y\to k} [\Gamma(y)].
\]
Indeed, it was shown that the homogeneous sea-state with autocorrelation $\Gamma(x-y)$ is unstable if, for some $X\in \mathbb{R},$ there exists $\Omega(X) \in \mathbb{C}$ so that
\beq\label{eq:Penr1}
1+\omega_0 k_0^2 \int\limits_k \frac{S(k+\frac{X}2)-S(k-\frac{X}2)}{\Omega + \frac{\omega_0}{4 k_0^2} kX} dk =0.
\ee
This was called an ``eigenvalue relation'' in \cite{Alber1978}; we will call it an ``instability condition''. In \cite{Athanassoulis2018}  it was further shown that if the instability condition does not hold then linear stability follows. The instability condition \eqref{eq:Penr1} itself can be refined in two ways: one concerns a technical issue related to $X=0.$ However, the more serious one is that in \eqref{eq:Penr1} we are asked to guarantee the existence or not of solutions for a non-linear system of two equations (the real and imaginary part of eq. \eqref{eq:Penr1}) in three real unknowns ($X,$ $\real \Omega,$ $\imag \Omega$).
This is not straightforward in general, and historically it has been a challenge in using the Alber equation more widely \cite{Gramstad2017}. In \cite{Athanassoulis2018} this condition is reformulated so that a more constructive way to check it can be found: by dividing both sides of the fraction by $X$ and setting 
\[
X'=\frac{X}{k_0}, \quad \Omega'=-\frac{\Omega 4 k_0}{X \omega_0}, \quad k'=\frac{k}{k_0},
\]
eq. \eqref{eq:Penr1} becomes
\[
\frac{1}{4\pi} = \frac{1}{\pi} k_0^3 \int\limits_{k'}  \frac{ \frac{S\big( (k'+\frac{X'}2)k_0 \big) - S\big( (k'-\frac{X'}2)k_0 \big)}{X'} }{\Omega'-k'} dk'
\]
This becomes very simple if we recognize it as {the Hilbert transform} of {the divided difference} of a {rescaled spectrum}. So, denoting
\beq\label{eq:defrescspec}
P(k) := k_0^3 S(kk_0), \quad D_XP(k) = \left\{ \begin{array}{cl}
\frac{P(k+\frac{X}2)-P(k-\frac{X}2)}{X}, & X\neq 0\\
P'(k), & X=0.
\end{array}
\right., \quad
\mathbb{H}[u](x)=\frac{1}\pi \mbox{p.v.} \int\limits_{t\in\R} \frac{u(t)}{x-t}dt
\ee
and dropping the primes,
the condition for instability finally becomes
\beq \label{eq:PenrHilsyn}
\exists X\in \R : \quad  
\exists \Omega=\Omega(X)\in \C  : \quad
\mathbb{H}[D_X P](\Omega) = \frac{1}{4\pi}.
\ee
The benefit is that now the argument principle\footnote{For a holomorphic function $f$ defined on a closed domain $A\subseteq \C,$ $f:A\to \C,$ it follows that $z\in f(A)$ if and only if the curve $f(\partial A)$ is circumscribed around $z\in \C.$} can be used  to reformulate this to a constructive condition that we can directly check \cite{Penrose1960,Athanassoulis2017,Athanassoulis2018}. To that end we will also need to introduce the {signal transform},
\beq\label{eq:signaltransdform}
\mathbb{S}[u](x):= \mathbb{H}[u](x) - i u(x);
\ee
we are now ready to state the equivalent instability condition:

\begin{definition}[Penrose-Alber condition] \label{def:penralbcond} Consider a homogeneous sea state with autocorrelation function $\Gamma(x-y)$ and carrier wavenumber $k_0,$ and denote $S(k)=\mathcal{F}_{y\to k} [\Gamma(y)],$ $P(k) := k_0^3 S(kk_0).$
Then, the sea state is Penrose-Alber unstable if 
\[
d(\overline\Gamma,\frac{1}{4\pi}) =0,
\]
where 
\beq\label{eq:defcurveG}
\bac
\Gamma_X:= \left\{\mathbb{S}[D_X P(\cdot)](t), \,\, t \in \R\right\}\cup \{0\}, \quad \overset{\circ}{\Gamma}_X=\{ z\in\C |  z \mbox{ enclosed by } \Gamma_X\}, \quad  
\overline\Gamma := \overline{\bigcup\limits_{X\in\R} \overset{\circ}{\Gamma}_X}.
\ea
\ee   
Any $X$ for which $1/4\pi \in \overline{\overset{\circ}{\Gamma}}_X$ is called an {unstable wavenumber}.
\end{definition}
This notion of instability is equivalent to condition \eqref{eq:PenrHilsyn} \cite{Athanassoulis2018}. The closed curves $\Gamma_X$ are vizualized for a concrete example in the top of Figure \ref{fig:1}.

%
%
%
%

\begin{figure}[!h]
\begin{center}
\includegraphics[width=0.96\textwidth]{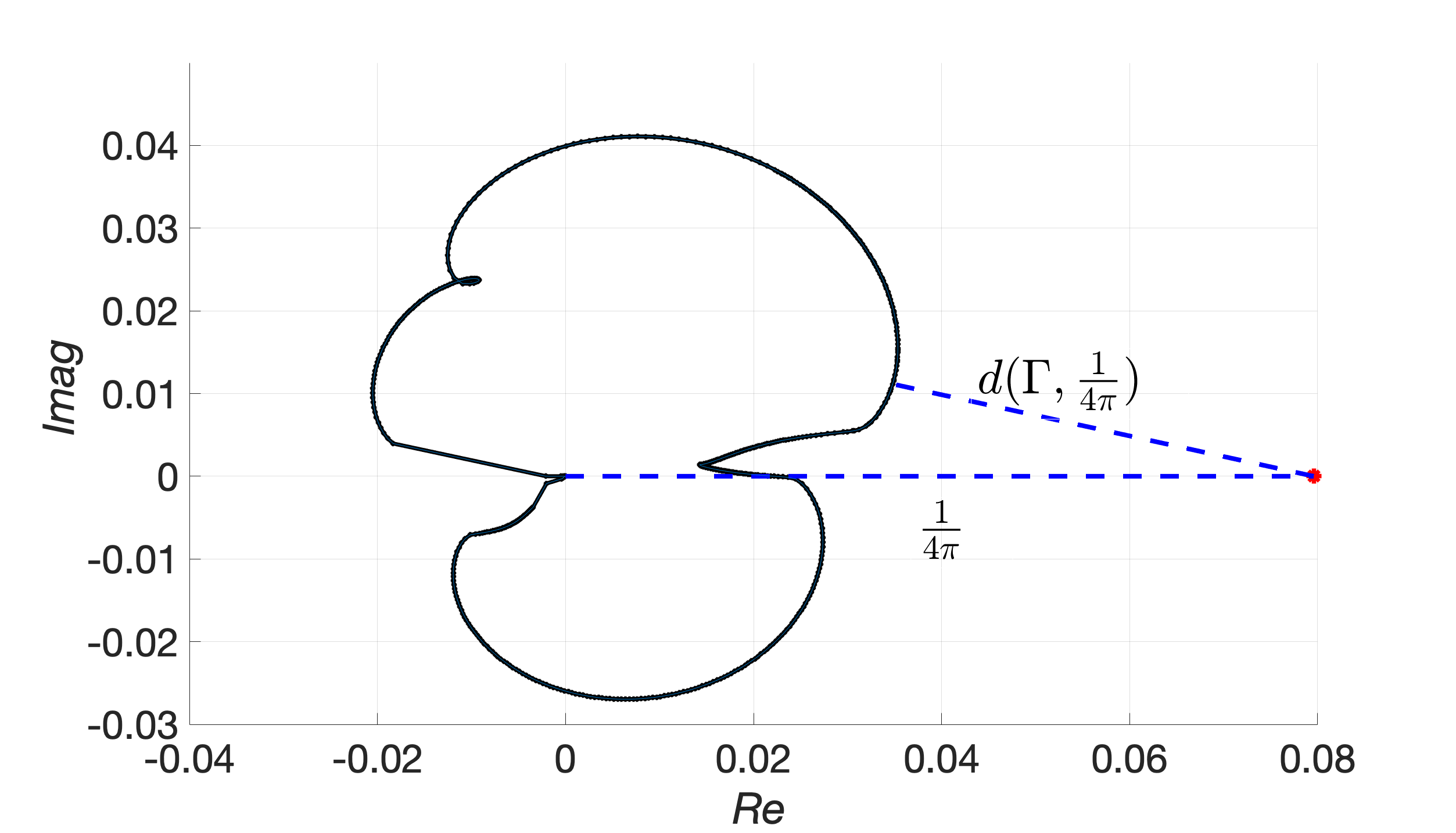} 
\caption{ 
We can quantify how close each spectrum is to being modulationally unstable by measuring the distance by
between $\Gamma_X$  and  $1/4\pi$ (plotted here as a red star). In practice, it suffices to do this for $X\approx 5\cdot 10^{-4},$ as much smaller values of $X$ yield similar results, and large values of $X$ lead to smaller and smaller curves, cf. Figure \ref{fig:1}.  So now we can say that spectrum no. 11300  has come something like $42\%$ of how far out it would need to be in order to be modulationally unstable. 
 }\label{fig:1b}
\end{center}
\end{figure}

%
%
%

We call this the Penrose-Alber instability condition after Alber's ``eigenvalue relation'' \eqref{eq:Penr1} \cite{Alber1978} and Penrose's introduction of the argument principle in an analogous problem in plasma \cite{Penrose1960}.
The point is that this formulation of the condition boils down to drawing closed curves on the complex plane, and looking whether the point $1/4\pi$ is inside or outside of them, cf. Figure  \ref{fig:1b}. It was used in \cite{Athanassoulis2018} to study parametric JONSWAP spectra; here we will apply it to non-parametric spectra as well.


%

\medskip

This stability-of-homogeneity question that was first studied in \cite{Alber1978} is in fact a variant of a much better known stability question. The {modulation instability (MI)} is well understood and widely studied  in a phase-resolved setting, for a plane wave background \cite{Zakharov2009}. It can easily be seen that a sequence of spectra approaching a plane wave $S(k)=a\delta(k-k_0)$ would become unstable in the Penrose-Alber sense. In fact, the Penrose-Alber instability is just the modulation instability for more general backgrounds than plane waves. A modulationally unstable sea state supports the rapid concentrations of energy -- a possible mechanism for the formation of rogue waves \cite{Onorato2009,Onorato2006,Onorato2013,Athanassoulis2017}.


Unlike the classical MI, where every plane-wave solution is always unstable, a spectrum can be stable or unstable in the Penrose-Alber sense. In the case of stability, the homogeneity of the sea-state is robust, and small perturbations will merely disperse. {In these cases, despite having a focusing nonlinearity and infinite energy present,  the dynamics are going to be dominated by the linear dispersion for all times.} This is a familiar leading order approximation in ocean waves, but it doesn't necessarily have a name to itself in an oceanographic context. As the rigorous stability analysis of \cite{Athanassoulis2018} highlighted, mathematically this stable regime looks exactly like what is called {Landau damping} for Vlasov equations (indeed this parallel was also drawn before \cite{Onorato2003}).

\subsection{Implications}

\subsubsection{Quantifying stability}

As was seen just above, a key feature of the Alber equation is a classification of a given spectrum as either stable or unstable. A careful look at the asymptotics can offer more nuance. For example, consider two sea states: one with a barely stable spectrum $S(k)$, and the other with a slightly perturbed version, e.g. $(1+\varepsilon)S(k),$  so that it becomes barely unstable. These sea states will behave similarly on physically realistic timescales -- as one would intuitively expect. In particular, there is not a violent bifurcation from Landau damping to modulation instability, but rather a gradual transition \cite{Athanassoulis2018,Athanassoulis2017}. On the other hand, as we will see in some detail here, a spectrum exhibiting Landau damping can be ``more stable'' than another  spectrum which also exhibits Landau damping. In short, the {effective stability} of a spectrum is better thought of as belonging to a continuous range of values rather than a binary ``stable / unstable'' classification. One of the results of this paper is a non-dimensional index quantifying this ``effective stability'', cf. Figures \ref{fig:1b}, \ref{fig:8}. 

This is important because, even in the presence of Landau damping (i.e. for ``stable'' spectra)  {only small enough inhomogeneities are guaranteed to disperse} if the nonlinearity is taken into account. Large inhomogeneities are not well understood, and could very well behave differently. This is a harder mathematical problem, because now by definition asymptotics are not sufficient. To better understand this question  there is a numerical requirement  (efficient and reliable solvers for the Alber equation both in the stable and unstable regime) as well as a data requirement (quantify how large inhomogeneities are in the ocean; this may be possible e.g. with X-band radar imaging \cite{Borge2004}). It is quite possible that, if a spectrum is close enough to MI  the likelihood  of nonlinear events under realistic perturbations may substantially increase, even if technically the spectrum still exhibits Landau damping.

\subsubsection{Nonparametric spectra}

The fundamental scaling of the problem shows that the vast majority of plausible sea states would be stable, in accordance with the well-known fact that linearized dynamics very often do a  good job in describing phase-averaged energy propagation. On the other hand, it seems that instability is  within reach, i.e. the scaling of the problem does not make instability  so far removed as to be considered impossible. This much was established in \cite{Athanassoulis2018,Gramstad2017,Ribal2013a}, working in the context of fitted parametric JONSWAP spectra. While fitted spectra are widely used, realistic spectra from the field come in many different shapes and forms.
In this paper we work with a set of nonparametric unidirectional spectra. We investigate whether they are stable or unstable in the Penrose-Alber sense, and proceed to examine how well the ``proximity to instability'' (defined more precisely in equation \eqref{def:pti}) correlates with the probability of rogue waves appearing in the given sea state.   

The dataset, methodology and results of our investigation are described in detail in Section \ref{sec:measspectra}. 

\subsubsection{Emergence of coherent structures}

\begin{figure}[h!]
\begin{center}
\includegraphics[width=0.7\textwidth]{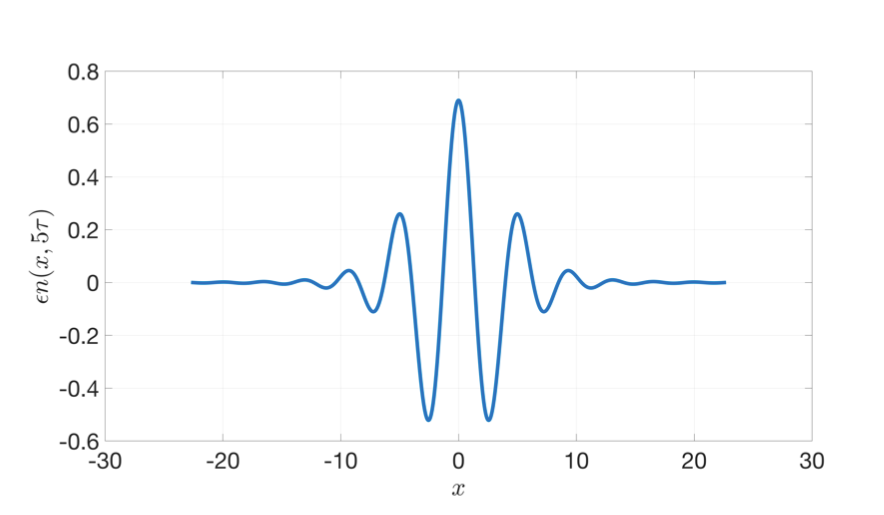} 
\caption{Predicted profile of emergent localized extreme events for NLS, computed according to the methodology of  \cite{Athanassoulis2017}. A virtually identical  universal profile of fully developed rogue waves is reported in Figure 2 of  \cite{Dematteis2017} for the largest extreme waves. 
This brings to mind the ``three sisters'' discussion \cite{Bitner-Gregersen2020,Magnusson2019,Muller2005,Nikolkina2011,Clauss2002,Ankiewicz2011}, or the Greek  $\tau\rho\iota\kappa\upsilon\mu\iota\alpha.$}
\label{fig:2}
\end{center}
\end{figure}

Another area with many open questions is what happens when instability arises. For modulationally unstable spectra with a small inhomogeneity, formal asymptotics can be used to describe the early evolution of the instability. A particular coherent structure emerges, determined by the unstable wavenumbers for the particular spectrum and their rate of growth. In that sense the coherent structure, at least in its early stages, is determined by the spectrum (and not by the inhomogeneity) and  the Penrose-Alber instability analysis suffices to predict it \cite{Athanassoulis2017}. Rather surprisingly, the same kind of universal coherent structure was reported in a fully numerical study by van den Eijden et al. \cite{Dematteis2017}, for the {fully  nonlinear stage of the instability}, cf. Figure \ref{fig:2}. This is a direction  with many open questions, where more work is needed.

\section{Main results} \label{sec:mainresults}

\subsection{The Alber equation for crossing seas}

The original derivation of the Alber equation \cite{Alber1978} allows for an oblique but small inhomogeneity on a unidirectional sea state. This does include some two dimensional aspects, however it does not allow for two large, different wavesystems, with different main directions of propagation, crossing. Such a situation is called ``crossing seas'', and
is recently attracting a lot of interest. In particular, it is thought that   modulation instability and rogue waves may be more prominent in crossing seas, but this is  still very far from fully understood
\cite{Gramstad2018,Steer2019,Hammack2005,Onorato2006,Onorato2010}. 

One way to study crossing seas is by deriving a coupled system of equations, each governing the evolution of one (quasi-uni-directional) wavefield. Let us denote the direction of propagation for the wavetrain $A$ by $\mathbf{k^A}=(k^A_1,k^A_2)$ and for the wavetrain $B$ by $\mathbf{k^B}=(k^B_1,k^B_2).$  The corresponding frequencies are 
\beq
\omega^A = \sqrt{g |\mathbf{k^A}|}, \qquad 
\omega^B = \sqrt{g |\mathbf{k^B}|}
\ee
for ocean waves (i.e. in the limit of infinite depth).
In \cite{Hammack2005} the following system of NLS equations for the envelopes of two wavetrains, $v^A,v^B,$ in two spatial dimensions is derived: 
\beq \label{eq:syst2nls_10}
i\frac{d}{dt} v^A+i\mathbf{C^{A}} \cdot \nabla_x v^A+ \alpha_1 \partial_{x_1}^{2}  v^A 
+ \beta_1 \partial_{x_2}^{2} v^A+ \gamma_1 \partial_{x_1}\partial_{x_2} v^A+\left(\xi_1|v^A|^{2}+ \zeta_1|v^B|^{2}\right) v^A=0,
\ee
\beq \label{eq:syst2nls_20}
i\frac{d}{dt} v^B+i\mathbf{C^{B}} \cdot \nabla_x v^B+ \alpha_2 \partial_{x_1}^{2}  v^B 
+ \beta_2 \partial_{x_2}^{2} v^B+ \gamma_2 \partial_{x_1}\partial_{x_2} v^B+\left(\xi_2|v^B|^{2}+ \zeta_2|v^A|^{2}\right) v^B=0.
\ee 
All the coefficients are completely determined in terms of $\mathbf{k^A},\mathbf{k^B}$ \cite{Hammack2005}.

\begin{theorem}\label{thrm:main} Consider the two crossing wavefileds of equations \eqref{eq:syst2nls_10}, \eqref{eq:syst2nls_20}, and assume moreover that their autocorrelations at $t=0$ can be written as
\beq\label{eq:assumeweakhiyg}
\E[v_A(\mathbf{x},t) \overline{v_A(\mathbf{y},t)}] = G^A(\mathbf{x}-\mathbf{y})+ \epsilon r_0(\mathbf{x},\mathbf{y}), \qquad 
\E[v_B(\mathbf{x},t) \overline{v_B(\mathbf{y},t)}] = G^B(\mathbf{x}-\mathbf{y})+ \epsilon s_0(\mathbf{x},\mathbf{y})
\ee
for some $\epsilon=o(1).$ 
Then the system \eqref{eq:syst2nls_10}, \eqref{eq:syst2nls_20} exhibits {modulation instability} if
\beq
\exists \mathbf{P} \in \R^2, \,\, \omega \in \C \,\,: \,\, (1-\xi_1 h^A(\mathbf{P},\omega))(1-\xi_2 h^B(\mathbf{P},\omega)) = \zeta_1\zeta_2 h^A(\mathbf{P},\omega)h^B(\mathbf{P},\omega)
\ee
where $h^A,h^B,\widetilde{n}_A^0,\widetilde{n}_B^0$ are defined in terms of the data of the problem in eq. \eqref{eq:techdefs}\footnote{Roughly speaking, $h^A$ and $h^B$ are transfer functions generated by the homogeneous backgrounds $\Gamma^A,$ $\Gamma^B,$ and $\widetilde{n}_A^0,\widetilde{n}_B^0$ express the free space (i.e. $\xi_1=\xi_2=\zeta_1=\zeta_2=0$) evolution of the initial inhomogeneities.}. In that case, inhomogeneities are expected to grow in time.

\bigskip

On the other hand, if
\beq\label{eq:stabilitycond8760}
\inf\limits_{\substack {\real \omega >0 \\ \mathbf{P}\in \mathbb{R}^2}} \Big|(1-\xi_1 h^A(\mathbf{P},\omega))(1-\xi_2 h^B(\mathbf{P},\omega)) - \zeta_1\zeta_2 h^A(\mathbf{P},\omega)h^B(\mathbf{P},\omega) \Big| =\kappa> 0
\ee
and
\beq\label{eq:stabilitycond8765}
\sup\limits_{\substack {\real \omega >0 \\ \mathbf{P}\in \mathbb{R}^2}} \big( |h^A(\mathbf{P},\omega)|+|h^B(\mathbf{P},\omega)| \big)  < +\infty
\ee
hold, then formally the problem exhibits {linear Landau damping}, i.e. inhomogeneities are expected to disperse and thus deviation from homogeneity to not grow noticeably.

\end{theorem}

The proof is found in Section \ref{sec:crossseas}.

\begin{remark}
Condition \eqref{eq:stabilitycond8760} is a Penrose-Alber condition for a system, and it is consistent with the condition for the scalar case. To see this assume that $v_B \to 0;$ then, according to the definition of $h^B$ (cf. eq. \eqref{eq:techdefs}) $h^B\to 0$ as well, and the condition becomes
\[
\inf\limits_{\substack {\real \omega >0 \\ \mathbf{P}\in \mathbb{R}^2}} \Big|1-\xi_1 h^A(\mathbf{P},\omega)\Big| =\kappa> 0
\]
which is exactly of the same type as requiring
\[
\inf\limits_{\substack {\real \Omega >0 \\ X\in \mathbb{R}}}
\Big|1 - 4\pi \mathbb{H}[D_X P](\Omega) \Big|=\kappa>0
\]
in the scalar case.

Just like in the scalar case, Condition \eqref{eq:stabilitycond8760}
 can also be resolved with the argument principle since,  for each $\mathbf{P}\in\mathbb{R}^2,$ it asks whether a holomorphic function attains the value $0$ on the right half-plane, $\inf_{{\real \omega >0, \, \mathbf{P}\in \mathbb{R}^2}} |F_\mathbf{P}(\omega)|>\kappa$.
Condition \eqref{eq:stabilitycond8765}
apparently boils down to general regularity conditions for the spectra $\widehat{\Gamma^A},$ $\widehat{\Gamma^B}.$ Working these ideas out in full detail will require the extension of some technical results to a non-standard ``two-dimensional Hilbert transform'' which arises here.
\end{remark}


\subsection{Stability of unidirectional non-parametric spectra and Proximity to Instability (PTI)} \label{sec:measspectra}

\subsubsection{The data}

For the present analysis we have used 100 realistic nonparametric spectra taken from the Norwegian Meteorological Institute's operational spectral wave model \cite{ReiEtAl2011JGR}, which is a third-generation wave model based on WAM. The model provides wave spectra every hour for many locations in the North Atlantic. For this study 100 spectra from one specific location south-west of the Norwegian coast ($56\degree\, 36'\,\text{N},\; 3\degree\, 12'\,\text{E}$) were used. The selected spectra consisted of 80 spectra randomly selected from the full database (26\,255 spectra covering the period from October 2016 to September 2019), as well as the 10 spectra having the largest mean wave steepness $\epsilon=H_sk_0/2$ and the 10 spectra with the largest BFI values. Recall that the Benjamin-Feir Index (BFI) was defined from the frequency spectrum $E(\omega)$ as
\[
  \BFI=\frac{\epsilon}{\sqrt{2}\delta_{\omega}},
\]
where $\delta_{\omega}$ is a measure for the frequency-bandwidth here defined in terms of Godas' peakedness factor $Q_p$, as suggested e.g. in \cite{JanBid2009REPORT}:
\[
  \delta_{\omega}=\frac{1}{Q_p\sqrt{\pi}} \quad \mbox{where} \quad Q_p=\frac{2}{m_0^2}\int \omega E^2(\omega)d\omega,
\]
and where $m_0=\int E(\omega)d\omega$ is the total energy of the spectrum.

Note that for the following analysis the original frequency spectra $E(\omega)$ were converted into wavenumber spectra $S(k)$ using the linear dispersion relation $\omega=\sqrt{gk}$.

\subsubsection{The algorithm for checking the instability condition}

Here we will implement the stability criterion of Definition \ref{def:penralbcond} to a number of non-parametric ocean spectra. First of all we will rescale them $S(k) \mapsto P(k):=k_0^3S(k\cdot k_0).$ Actually the selection of $k_0$ is a nontrivial question. For example, one could plausibly use the mode, the mean or the median wavenumber. For narrow, unimodal spectra the choice would make very little difference, but for more irregular shapes, including e.g. bimodal spectra, the resulting differences might be noticeable. In this case, we use the $k_0$ provided in the dataset by the Norwegian Meteorological Institute, and which has been used for the computation of the BFI, steepness contained in the dataset etc. 

Once $k_0$ is determined and the rescaling is complete, we will need to  interpolate the discrete data to a finer grid. (The original discrete spectra come sampled in 36 non-uniformly spaced wavenumbers.)  This will be crucial in using effective quadrature methods. Splines are well suited to this task; the additional requirement is to minimize overshooting at maxima, as this could make a big difference with regard to our investigation. To that end we use  \texttt{pchip}, a piecewise polynomial interpolation routine in MATLAB that minimizes overshoot. 


The non-trivial part of the computation is the  Hilbert transform, which is a singular integral. To that end we will use the Sokhotski-Plemelj formula,
\[
\forall u \in C(\R) \cap L^1(\R) \qquad 
\lim\limits_{\eta\to 0^+}\frac{1}{\pi} \int\limits_{s} \frac{u(s)}{t-s-(\sqrt{-1}) \eta }ds = \mathbb{H}[u](t)-i u(t)=  \mathbb{S}[u](t)
\]
and truncate the limit by taking an appropriate $\eta=\texttt{compl\_tol}\ll 1.$ After extensive testing it is found that the result doesn't change noticeably once $\eta\approx 10^{-4}$ That way we avoid the singular integral. A detailed pseudocode for how the curve $\Gamma_X$ is generated is presented in Algorithm \ref{alg:1}. 

It is a moderately heavy computation if a good approximation for all of $\Gamma_X$ is required, as a few million points are typically required in order to achieve stringent error tolerances ($\sim 10^{-2}$ relative error tolerance or $\sim 10^{-6}$ absolute error tolerance). However, one can check stability more quickly, by checking only the points $t_*$ where $D_X P(t_*)=0;$ these are the points where the curve $\Gamma_X$ crosses the real axis. If it never crosses the real axis to the right of $1/4\pi,$  then topologically $1/4\pi$ cannot possibly be in the interior of the curve.

Some more details involve the selection of $X;$ it is found that $X\ll 1$ will produce the curves that have the most chance of coming closer to $1/4\pi,$ as when $X$ increases the curves $\Gamma_X$ shrink to zero, see Figure \ref{fig:1}. Numerical testing shows that $X\approx 10^{-4}$ gives a good picture of what happens as $X\to 0.$ So, if for $X\approx 10^{-4}\,$ $\Gamma_X$ is not winding around $1/4\pi,$ nor coming very close to it, we can accept the spectrum as stable. 

%
%

\begin{algorithm}
\begin{itemize}
\item[ ]	 {\bf Input } $S_j,$ $k_j$   {\em (Sampled values of wavenumber-resolved spectrum, $k_j \in [0.00479,3.78].$)} 
\item[ ] {\bf Rescale} $(S_j,k_j) \mapsto (k_0^3 S_j, k_j/k_0)=(P_j,\xi_j)$
 {\em  ($k_0$ is simply taken to be the peak wavenumber.)} 
\item[ ]  {\bf Interpolate } $P(\xi)$ 
\item[ ] {\bf Set} \quad \texttt{compl\_tol}$\sim 10^{-4}$, \,\,  \texttt{rel\_tol}$\sim10^{-2}$, \,\,  \texttt{abs\_tol}$\sim10^{-6}$
\item[ ] {\bf For} \quad $t_i:=t_{min}$ \,\, {\bf to}\,\, $t_{max}$\,\, {\bf step}\,\, $\delta t$ 
\begin{itemize}
\item[ ] {\bf While} \texttt{rel\_err} $>$ \texttt{rel\_tol} {\bf AND} \texttt{abs\_err} $>$ \texttt{abs\_tol}
\begin{itemize}
\item[ ] {\bf Integrate} 
\[
I= \frac{1}{\pi} \int\limits_{s} \frac{\frac{P(s+\frac{X}2)-P(s-\frac{X}2)}{X}}{t_i-s-(\sqrt{-1}) \, \mbox{\texttt{compl\_tol}} }ds
\]
using composite Simpson on two quadrature grids: a finer one and a coarser one, generating two approximations, \texttt{I\_fine} and \texttt{I\_coarse}. {\em (Fine grid has $3\times$the number of points compared to coarse grid.) }
\item[ ] {\bf Set} $\texttt{rel\_err} = \frac{|\texttt{I\_fine}-\texttt{I\_coarse}|}{\texttt{I\_fine}},$ \,\,\, $\texttt{abs\_err} = |\texttt{I\_fine}-\texttt{I\_coarse}|.$ 
\end{itemize}
\item[ ] {\bf End While}
\item[ ] {\bf Set} $\Gamma(t_i)=\texttt{I\_fine}$
\end{itemize}
\item[ ] {\bf End For}
\item[ ] {\bf Plot} the line $\Big(\real \Gamma_X(t_i), \imag \Gamma_X(t_i) \Big),$ $i=1,2,...$ and the point $(\frac{1}{4\pi},0)$ 
\item[ ] {\bf Check} whether $(\frac{1}{4\pi},0)$ is inside $\Gamma_X$
 {\em (This can be done with the MATLAB \texttt{inpolygon} function.)}
\end{itemize}
\caption{ Pseudo-code for the computation of $\Gamma_X$ }\label{alg:1}
\end{algorithm}

\subsubsection{Summary of the results} \label{sec:summary}

Most of the non-parametric spectra examined were found to be stable, i.e. exhibit Landau damping. However, three spectra were found to be modulationally unstable, and a handful more were extremely close to being unstable. 

\begin{figure}[!h]
\begin{center}
\includegraphics[width=0.48\textwidth]{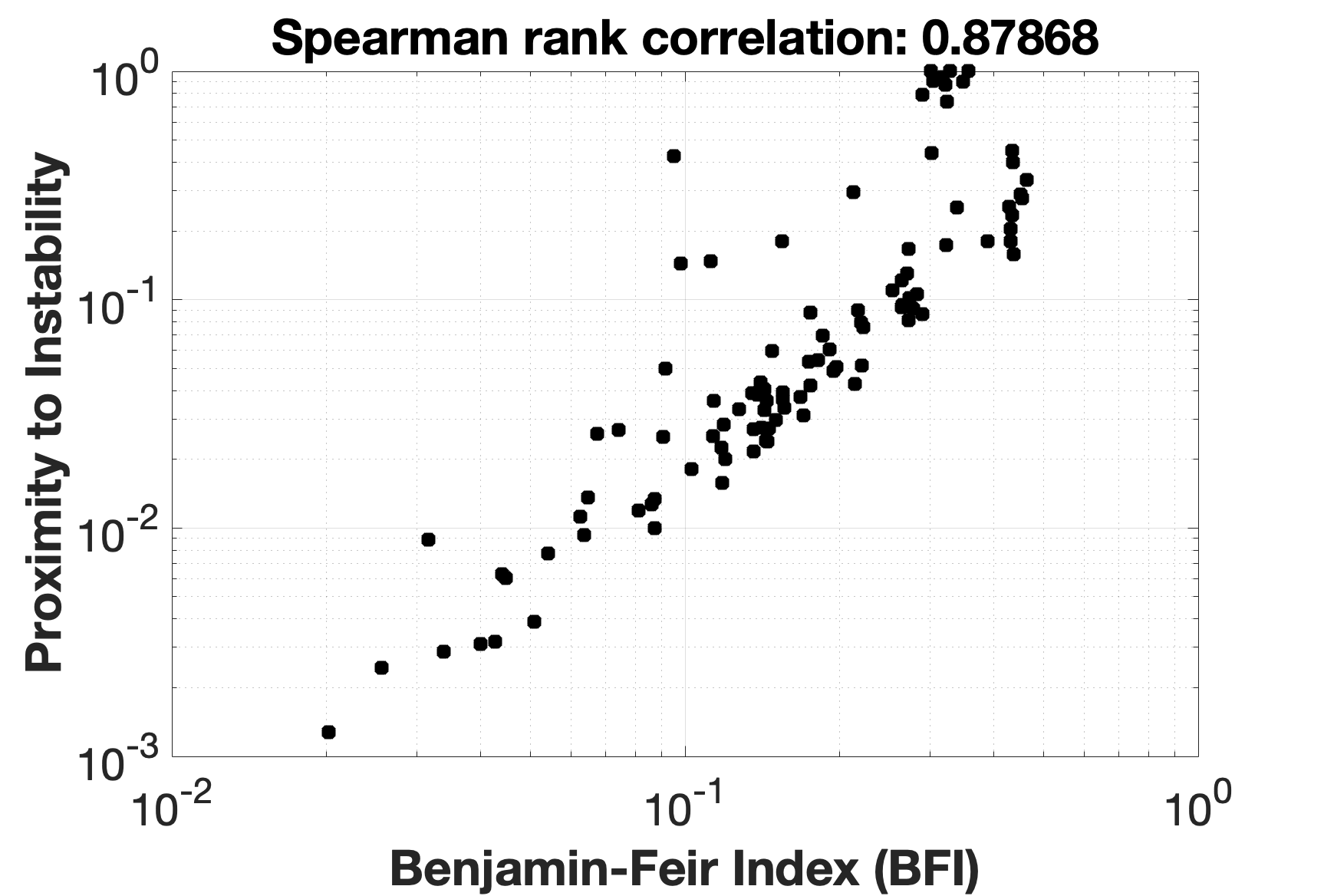} \,\,
\includegraphics[width=0.48\textwidth]{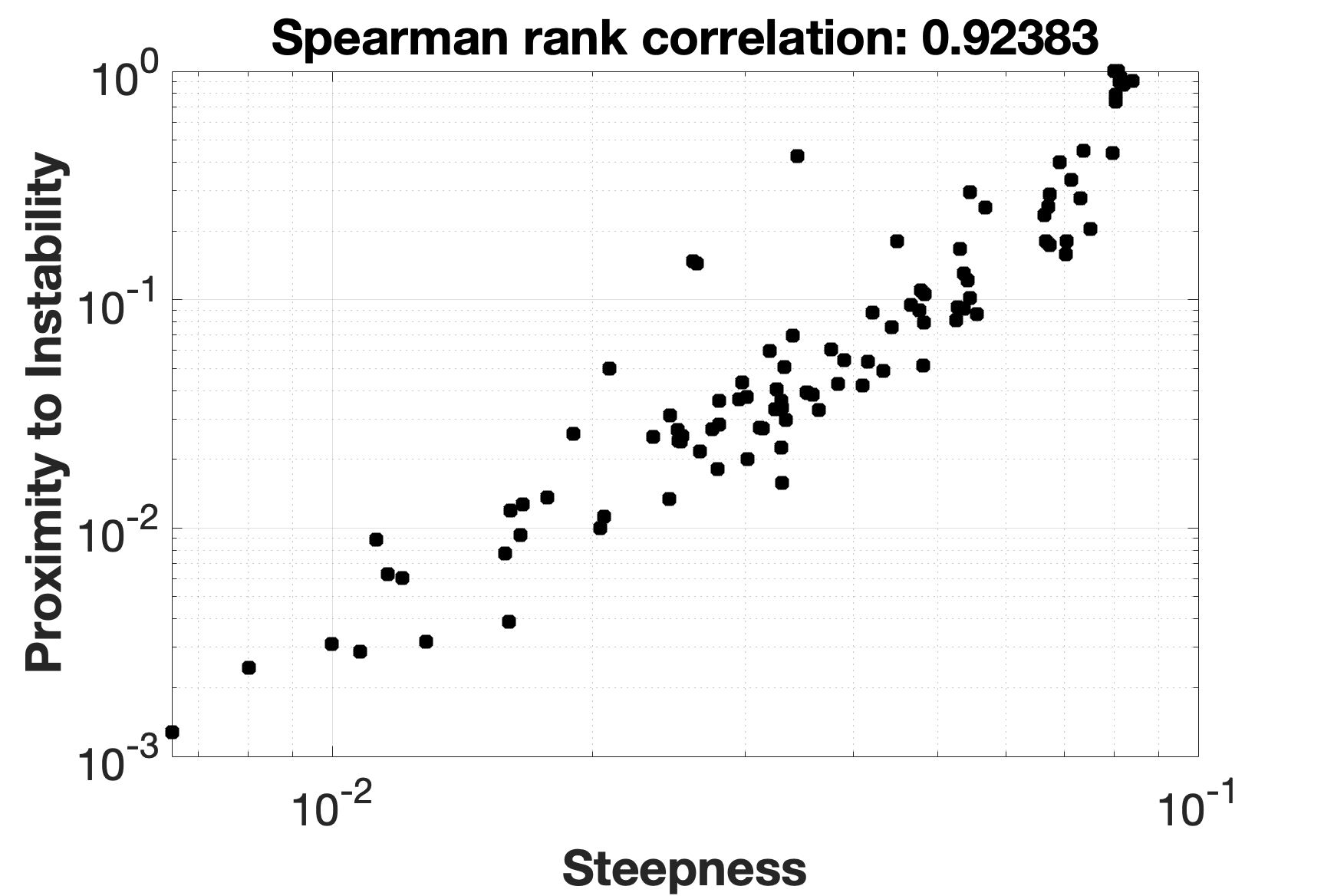} 
\caption{ The new metric of PTI offers a quantitative and nondimensional way to assess how far a spectrum is from being modulationally unstable.  The Benjamin-Feir Index (BFI), essentially a rescaled version of steepness, was introduced with a similar task in mind; so it is only natural to ask how well the two are correlated. In the {\bf  left} graph above we see a log-log scatterplot of the BFI v. PTI. In the {\bf  right} graph we see a similar scatterplot of a representative wave steepness, $\epsilon=H_sk_0/2,$  v. PTI.
 }\label{fig:8}
\end{center}
\end{figure}

\begin{figure}[!h]
\begin{center}
\includegraphics[width=0.48\textwidth]{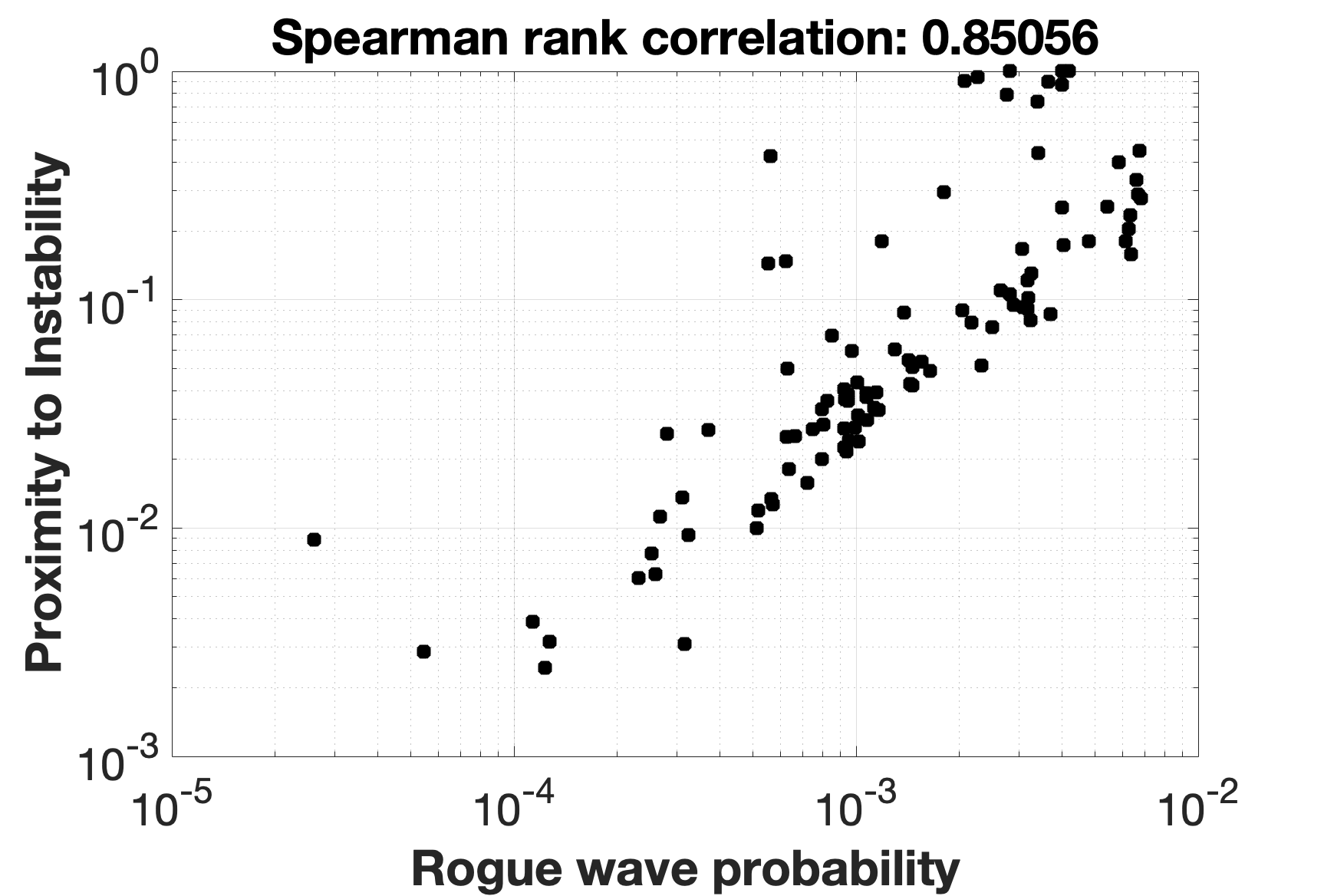} \,\,
\includegraphics[width=0.48\textwidth]{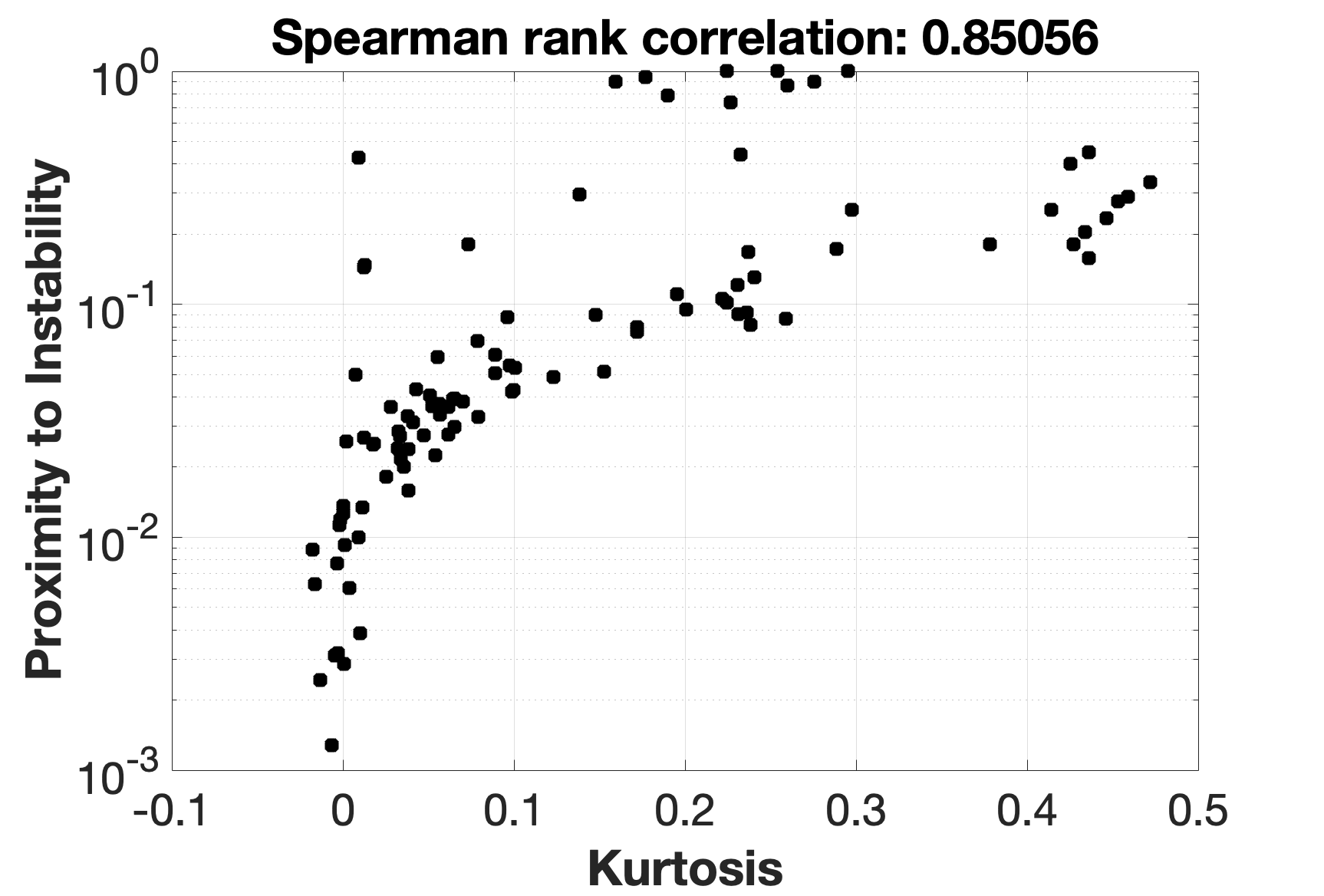} 
\caption{ Comparison of PTI to Monte Carlo results, please see Section \ref{sec:summary} for more details on how these are defined and computed.  {\bf Left} Rogue wave probability (understood as $P(\mbox{wave crest }>H_s)$) v. PTI.  {\bf Right} Kurtosis v. PTI. See Section \ref{sec:summary} for more details and context.
 }\label{fig:82}
\end{center}
\end{figure}

We also define as proximity to instability (PTI) the quantity
\beq\label{def:pti}
\mbox{PTI} = 1 - \frac{d(\overline{\Gamma},\frac{1}{4\pi})}{\frac{1}{4\pi}};
\ee 
this is 1 for any modulationally unstable spectrum, and $0$ for the zero spectrum. It provides a non-dimensional way to quantify how close a spectrum comes to being modulationally unstable in the sense of Definition \ref{def:penralbcond}. We will compare this with quantities of interest for the same spectrum (cf. Figure \ref{fig:8}), along with a Monte Carlo estimation of the likelihood of rogue waves and nonlinear events (cf. Figure \ref{fig:82})

Indeed, for each of the 100 spectra selected for this study, we have run numerical simulations with the Higher Order Spectral Method (HOSM) \cite{dommermuth1987high,west1987new} in a Monte-Carlo approach where each spectrum were simulated 100 times with different initial random phases and random amplitudes each run. That is, for a given spectrum $S(k)$ the initial surface elevation is in the form
\begin{align}
\eta(x,t=0)=\real\sum_{j=1}^{n}A_j\exp{(i k_jx)} \quad \mbox{where}\quad A_j=Z_j\sqrt{2S(k_j)\Delta k_j}
\end{align}
where $\Delta k_j=k_{j+1}-k_{j}$ is the grid spacing between the discrete wavenumbers and where $Z_j$ are independent complex standard normal variables. That is, the real and imaginary parts of $Z_j$ are independent normally distributed random variables with zero mean and variance $1/2$, meaning that $|Z_j|$ are Rayleigh distributed with parameter $\sigma=1/\sqrt{2}$ so that $E[|Z_j|^2]=1$ and $E[|A_j|^2]=2S(k_j)\Delta k_j$, and the phases $Arg(Z_j)$ are uniformly distributed on $[0, 2\pi)$. 

HOSM applies a regular discretization of the wavenumbers so that $\Delta k=2\pi/x_{max}=2\pi/n\Delta x$, leading to a periodic domain of length $x_{max}$ in space. In the present simulations we have used $n=1024$, representing wavenumbers up to $k_{max}=8k_0.$. This means that $\Delta k=2k_{max}/n$, corresponding to $\Delta x=\lambda_0/16$ and $x_{max}=64\lambda_p$, where $\lambda_0=2\pi/k_0$ is the reference wavelength.

Each simulation was run for 30 minutes, from which time series of surface elevation were extracted from four locations distributed over the simulation domain. Thus, 200 hours of surface elevation time series were obtained for each sea state. Since here we are interested in relations between the instability (or proximity of such) obtained from the stochastic approach (i.e. Alber equation) and the occurrence of rogue and extreme events in random realizations of the sea states (i.e. phase-resolved Monte-Carlo simulations), we will consider the following parameters related to the occurrence of extreme wave events: {\bf the sea surface kurtosis} and {\bf the probability of extreme wave crests}.

{The sea surface kurtosis} is a measure for how much the tail of the distribution deviates from Gaussian statistics. For a Gaussian distribution the kurtosis is equal to zero, while a positive value for kurtosis  indicates more large events than in a Gaussian population. Hence, the kurtosis is often used as an indicator for the probability of extreme and rogue waves. 

Secondly we consider {the probability that a wave crest, defined as the maximum between each zero-crossing of the surface elevation, exceeds  the significant wave height ${H_s}$; ${P(C>H_s)}$\footnote{Recall that wave height is defined as the difference between a local maximum and a subsequent local minimum, i.e. in a sine wave it would be equal to two wave crests. Moreover, the significant wave height $H_s$ is defined as the mean value of the highest $33\%$ of wave heights.}. {This is a common definition of a rogue wave, which under linear \& gaussian assumptions (Rayleigh distributed crests) has the probability $P(C>H_s)=\exp{(-8)}$.} The probability of crest exceedance is estimated as the relative frequency of individual crests in the HOSM simulations that exceed the threshold $H_s$. For brevity we will refer to this probability simply as ``rogue wave probability''.

\section{Proof of Theorem \ref{thrm:main}}\label{sec:crossseas}

First of all observe that, by virtue of a simple gauge transform,
\beq\label{eq:gaugedef}
v^A(x_1,x_2,t) = e^{i (\kappa^A x_1 + \lambda^A x_2 + \tau^A t)} A(x_1,x_2,t), \quad 
v^B(x_1,x_2,t) = e^{i (\kappa^B x_1 + \lambda^B x_2 + \tau^B t)} B(x_1,x_2,t)
\ee
we can eliminate the  terms with the first derivatives\footnote{
By plugging \eqref{eq:gaugedef} in  eq. \eqref{eq:syst2nls_10} we see that eliminating the linear terms containing $A$ and $\nabla A$  leads to
\[
\begin{aligned}
\kappa^A 2\alpha_1+ \lambda^A\gamma_1 & =-C^A_1,\\
\kappa^A \gamma_1 + \lambda^A 2 \beta_1 & = -C^A_2,\\
\tau & =-C^A_1\kappa^A-C^A_2\lambda^A-\alpha_1(\kappa^A)^2-\beta_1(\lambda^A)^2-\gamma_1\kappa^A\lambda^A.
\end{aligned}
\] The first two lines is a $2x2$ system, and for the third, nonlinear equation we only need to substitute the obtained $\kappa^A,\lambda^A$ values from above. The same steps lead to a symmetric expression for $B.$}; without loss of generality we will work with the simplified system
\beq \label{eq:syst2nls_1}
i\frac{d}{dt} A+ \alpha_1 \partial_{x_1}^{2}  A 
+ \beta_1 \partial_{x_2}^{2} A+ \gamma_1 \partial_{x_1}\partial_{x_2} A+\left(\xi_1|A|^{2}+ \zeta_1|B|^{2}\right) A=0,
\ee
\beq \label{eq:syst2nls_2}
i\frac{d}{dt} B+ \alpha_2 \partial_{x_1}^{2}  B 
+ \beta_2 \partial_{x_2}^{2} B+ \gamma_2 \partial_{x_1}\partial_{x_2} B+\left(\xi_2|B|^{2}+ \zeta_2|A|^{2}\right) B=0.
\ee

Special symmetric cases of this system which simplify the coefficients have been studied  in \cite{Onorato2006,Steer2019}, but  we can in fact derive an Alber system for the general case.

To proceed we will use some shorthand notations, namely
\beq
\bac
\mathbf{x}=(x_1,x_2), \quad \mathbf{y}=(y_1,y_2), \quad  

L^j_{\mathbf{x}} = \alpha_j \partial_{x_1}^{2} 
+ \beta_j \partial_{x_2}^{2}+ \gamma_j \partial_{x_1 x_2}, \quad j \in \{ 1,2\}.
\ea
\ee
Now equation \eqref{eq:syst2nls_1} can be compactified as
\[
i\frac{\partial A}{\partial t}+ L^1_{\mathbf{x}}A+\left(\xi_1|A|^{2}+ \zeta_1|B|^{2}\right) A=0 
\]
and similarly for eq. \eqref{eq:syst2nls_2} with $L^2_{\mathbf{x}}$.

At this point let us introduce the autocorrelation functions for $A,B;$ 
\beq
R_{A}(\mathbf{x},\mathbf{y},t) = \E[A(\mathbf{x},t)\overline{A}(\mathbf{y},t)], \qquad R_{B}(\mathbf{x},\mathbf{y},t) = \E[B(\mathbf{x},t)\overline{B}(\mathbf{y},t)].
\ee
Observe that the relationship with the autocorrelations of the original $v^A,$ $v^B$ variables is
\beq\label{eq:convsm}
R_{A}(\mathbf{x},\mathbf{y},t) =e^{-i(\kappa^A(x_1-y_1)+\lambda^A(x_2-y_2))}\E[v^A(\mathbf{x},t)\overline{v^A(\mathbf{y},t)}]
\ee
and similarly for $R_B.$

By a straightforward computation we see that
\[
\bac
\frac{d}{dt} R_A(\mathbf{x},\mathbf{y},t) = \E[A(\mathbf{x},t)\, \frac{d}{dt}\overline{A}(\mathbf{y},t) \,\, + \,\, \overline{A}(\mathbf{y},t)\, \frac{d}{dt}{A}(\mathbf{x},t)]=\\

{ } \\

= \E\left[{ A(\mathbf{x},t)
\left(
-i L^1_{\mathbf{y}}\overline{A}(\mathbf{y},t)-i\left(\xi_1|A(\mathbf{y},t)|^{2}+ \zeta_1|B(\mathbf{y},t)|^{2}\right) \overline{A}(\mathbf{y},t)
\right)
}\right] + \\

+ \E\left[{ \overline{A}(\mathbf{y},t)
\left(
 iL^1_{\mathbf{x}}A(\mathbf{x},t)+i\left(\xi_1|A(\mathbf{x},t)|^{2}+ \zeta_1|B(\mathbf{x},t)|^{2}\right) A(\mathbf{x},t)
\right)
}\right] =
\\ { } \\

= i (L^1_{\mathbf{x}} - L^1_{\mathbf{y}}) R_A(\mathbf{x},\mathbf{y},t) + \hfill \\

+ i \xi_1\E[|A(\mathbf{x},t)|^2 A(\mathbf{x},t)\overline{A}(\mathbf{y},t)]
- i \xi_1\E[|A(\mathbf{y},t)|^2 \overline{A}(\mathbf{y},t){A}(\mathbf{x},t)]+ \\

+ i\zeta_1 \E[|B(\mathbf{x},t)|^2A(\mathbf{x},t)\overline{A}(\mathbf{y},t)] 
- i\zeta_1 \E[|B(\mathbf{y},t)|^2\overline{A}(\mathbf{y},t){A}(\mathbf{x},t)].

\ea
\]
The last two lines consist of  fourth order stochastic moments, and it is these terms that will have to be approximated by some closure scheme. For the fourth order moments involving $A$ only, we will use the same idea as in the standard Alber equation; namely we will use the fact that for Gaussian processes the relationship
\beq
\E[|A(\mathbf{x},t)|^2 A(\mathbf{x},t)\overline{A}(\mathbf{y},t)] = 2 R_A(\mathbf{x},\mathbf{x},t) R_A(\mathbf{x},\mathbf{y},t)
\ee
holds (cf. Theorem \ref{thrm:ComplexIsserlis}). This now provides the Gaussian closure for all terms of the same form, since
\beq
\E[|A(\mathbf{y},t)|^2 A(\mathbf{x},t)\overline{A}(\mathbf{y},t)] = 2 R_A(\mathbf{y},\mathbf{y},t) R_A(\mathbf{x},\mathbf{y},t)
\ee
follows by the same argument. 
The new kinds of terms that come into play for the first time are the joint $A-B$ moments. In physical terms, $A$ and $B$ represent two different wavetrains meeting in the ocean, each having been generated and propagated independently from each other. It is thus a reasonable assumption that they are stochastically independent. In that case, the joint moments simplify to
\beq
\bac
\E[|B(\mathbf{x},t)|^2A(\mathbf{x},t)\overline{A}(\mathbf{y},t)] = R_B(\mathbf{x},\mathbf{x},t)R_A(\mathbf{x},\mathbf{y},t), \\ { } \\
\E[|B(\mathbf{y},t)|^2\overline{A}(\mathbf{y},t){A}(\mathbf{x},t)] = R_B(\mathbf{y},\mathbf{y},t)R_A(\mathbf{x},\mathbf{y},t).
\ea
\ee

So, using these closures we finally end up with the following equation (whenever the independent variables are not shown explicitly, it is understood that $R_A:=R_A(\mathbf{x},\mathbf{y},t)$):
\beq \label{eq:AlberSystem_1}
\bac
i\frac{d}{dt} R_A +
(L^1_{\mathbf{x}} - L^1_{\mathbf{y}}) R_A  + \hfill \\

+2\xi_1 \left( R_A(\mathbf{x},\mathbf{x},t) - R_A(\mathbf{y},\mathbf{y},t) \right) R_A + 2\zeta_1 \left( R_B(\mathbf{x},\mathbf{x},t) - R_B(\mathbf{y},\mathbf{y},t)  \right) R_A=0.
\ea
\ee
By the symmetry of equations \eqref{eq:syst2nls_1}, \eqref{eq:syst2nls_2} one readily sees that, using the same kinds of closures,
\beq \label{eq:AlberSystem_2}
\bac
i\frac{d}{dt} R_B +
(L^2_{\mathbf{x}} - L^2_{\mathbf{y}}) R_B  + \hfill \\

+2\xi_2 \left( R_B(\mathbf{x},\mathbf{x},t) - R_B(\mathbf{y},\mathbf{y},t) \right) R_B + 2\zeta_2 \left( R_A(\mathbf{x},\mathbf{x},t) - R_A(\mathbf{y},\mathbf{y},t)  \right) R_B=0.
\ea
\ee

Now, we invoke the assumption of eq. \eqref{eq:assumeweakhiyg}, i.e.
\beq\label{eqL:leadordstate}
R_A(\mathbf{x},\mathbf{y},t) = \Gamma^A(\mathbf{x}-\mathbf{y}) + \epsilon \rho(\mathbf{x},\mathbf{y},t), \quad 
R_B(\mathbf{x},\mathbf{y},t) = \Gamma^B(\mathbf{x}-\mathbf{y}) + \epsilon \sigma(\mathbf{x},\mathbf{y},t)
\ee
where now $\rho,\sigma$ are the inhomogeneous components  of the autocorrelation\footnote{Recall that the conversion between second moments of $A,B$ and $v^A,v^B$ follows eq. \eqref{eq:convsm}.} Observe that, under the condition \eqref{eqL:leadordstate}, the system \eqref{eq:AlberSystem_1}, \eqref{eq:AlberSystem_2} is equivalent to (like earlier, $\rho:=\rho(\mathbf{x},\mathbf{y},t)$ and similarly for $\sigma$)
\beq \label{eq:AlberSystem_11}
\bac
i\frac{d}{dt} \rho +
(L^1_{\mathbf{x}} - L^1_{\mathbf{y}}) \rho  + \hfill \\

+2\xi_1 \left( \rho(\mathbf{x},\mathbf{x},t) - \rho(\mathbf{y},\mathbf{y},t) \right) (\Gamma^A(\mathbf{x}-\mathbf{y}) + \epsilon \rho) + 2\zeta_1 \left( \sigma(\mathbf{x},\mathbf{x},t) - \sigma(\mathbf{y},\mathbf{y},t)  \right) (\Gamma^A(\mathbf{x}-\mathbf{y}) + \epsilon \rho)=0, \\ { } \\

i\frac{d}{dt} \sigma +
(L^2_{\mathbf{x}} - L^2_{\mathbf{y}}) \sigma  + \hfill \\

+2\xi_2 \left( \sigma(\mathbf{x},\mathbf{x},t) - \sigma(\mathbf{y},\mathbf{y},t) \right) (\Gamma^B(\mathbf{x}-\mathbf{y}) + \epsilon \sigma) + 2\zeta_2 \left( \rho(\mathbf{x},\mathbf{x},t) - \rho(\mathbf{y},\mathbf{y},t)  \right)  (\Gamma^B(\mathbf{x}-\mathbf{y}) + \epsilon \sigma)=0,
\ea
\ee
and $\rho=\sigma=0$ is a solution. Following the terminology of Wigner equations, the terms 
\beq
n_A(\mathbf{y},t):=\rho(\mathbf{y},\mathbf{y},t), \qquad 
n_B(\mathbf{y},t):=\sigma(\mathbf{y},\mathbf{y},t)
\ee
will be called {position densities}; they are real valued, and they control the RMS amplitude of the inhomogeneities, for each wavetrain, on the coordinates $(\mathbf{y},t).$

The question that we will focus on is whether $\rho=\sigma=0$ is a linearly stable solution. To that end we will consider non-zero initial data $\rho(\mathbf{x},\mathbf{y},0), \sigma(\mathbf{x},\mathbf{y},0).$ 
We will also need to introduce the change of variables
\beq
\bac
\rho(\mathbf{x},\mathbf{y},t) = \widecheck{f}(\mathbf{p},\mathbf{q},t), \qquad 
\sigma(\mathbf{x},\mathbf{y},t) = \widecheck{g}(\mathbf{p},\mathbf{q},t) \\ { } \\

\mathbf{p} = \frac{\mathbf{x}+\mathbf{y}}2, \qquad \mathbf{q} = \mathbf{x}-\mathbf{y}
\ea
\ee
which leads to
\beq
\bac
\nabla_\mathbf{x} = \frac{1}2\nabla_\mathbf{p} + \nabla_\mathbf{q}, \qquad 
\nabla_\mathbf{y} = \frac{1}2\nabla_\mathbf{p} - \nabla_\mathbf{q}, \\ { } \\

L^j_\mathbf{x} 
=\alpha_j (\frac{1}2\partial_{p_1} + \partial_{q_1})^2 + \beta_j(\frac{1}2\partial_{p_2} + \partial_{q_2})^2 + \gamma_j (\frac{1}2\partial_{p_1} + \partial_{q_1}) (\frac{1}2\partial_{p_2} + \partial_{q_2}), \\

L^j_\mathbf{y} 
=\alpha_j (\frac{1}2\partial_{p_1} - \partial_{q_1})^2 + \beta_j(\frac{1}2\partial_{p_2} - \partial_{q_2})^2 + \gamma_j (\frac{1}2\partial_{p_1} - \partial_{q_1}) (\frac{1}2\partial_{p_2} - \partial_{q_2}), \\

\\ { } \\

\mathbf{x}=\mathbf{p}+\frac{\mathbf{q}}2, \qquad 
\mathbf{y}=\mathbf{p}-\frac{\mathbf{q}}2, \\ { } \\

n_A(\mathbf{x},t) = \rho(\mathbf{x},\mathbf{x},t)=\widecheck{f}(\mathbf{x},0,t)=\widecheck{f}(\mathbf{p}+\frac{\mathbf{q}}2,0,t), \qquad n_A(\mathbf{y},t)=\rho(\mathbf{y},\mathbf{y},t) = \widecheck{f}(\mathbf{p}-\frac{\mathbf{q}}2,0,t),
\ea
\ee
and similarly for $n_B.$ So
finally we arrive at the following system in the $\mathbf{p},\mathbf{q}$ variables ($\widecheck{f}:=\widecheck{f}(\mathbf{p},\mathbf{q},t),$ $\widecheck{g}:=\widecheck{g}(\mathbf{p},\mathbf{q},t)$):
\beq
\bac
i \frac{d}{dt} \widecheck{f} + \left[ 2\alpha_1 \partial_{p_1} \partial_{q_1} + 2\beta_1 \partial_{p_2}\partial_{q_2} + \gamma_1(\partial_{p_1}\partial_{q_2} + \partial_{p_2}\partial_{q_1}) \right] \widecheck{f} + \hfill \\

+ 2\left[
\xi_1\big(\widecheck{f}(\mathbf{p}+\frac{\mathbf{q}}2,0,t)-\widecheck{f}(\mathbf{p}-\frac{\mathbf{q}}2,0,t)\big) + \zeta_1
\big(\widecheck{g}(\mathbf{p}+\frac{\mathbf{q}}2,0,t)-\widecheck{g}(\mathbf{p}-\frac{\mathbf{q}}2,0,t)\big) \right] (\Gamma^A(\mathbf{q})+\epsilon \widecheck{f}) = 0, \\ { } \\

i \frac{d}{dt} \widecheck{g} + \left[ 2\alpha_2 \partial_{p_1} \partial_{q_1} + 2\beta_2 \partial_{p_2}\partial_{q_2} + \gamma_2(\partial_{p_1}\partial_{q_2} + \partial_{p_2}\partial_{q_1}) \right] \widecheck{g} + \hfill \\

+ 2\left[
\xi_2\big(\widecheck{g}(\mathbf{p}+\frac{\mathbf{q}}2,0,t)-\widecheck{g}(\mathbf{p}-\frac{\mathbf{q}}2,0,t)\big) + \zeta_2
\big(\widecheck{f}(\mathbf{p}+\frac{\mathbf{q}}2,0,t)-\widecheck{f}(\mathbf{p}-\frac{\mathbf{q}}2,0,t)\big) \right] (\Gamma^B(\mathbf{q})+\epsilon \widecheck{g}) = 0
\ea
\ee
Finally by and linearising (setting $\epsilon=0$) and taking the Fourier transform in both variables, 
\beq
f(\mathbf{P},\mathbf{Q},t) = \mathcal{F}_{\mathbf{p},\mathbf{q} \to \mathbf{P},\mathbf{Q}}[ \widecheck{f}(\mathbf{p},\mathbf{q})], \qquad 
g(\mathbf{P},\mathbf{Q},t) = \mathcal{F}_{\mathbf{p},\mathbf{q} \to \mathbf{P},\mathbf{Q}} [\widecheck{g}(\mathbf{p},\mathbf{q})],
\ee
we obtain
\beq
\bac
i \frac{d}{dt} {f} -4\pi^2 \left[ 2\alpha_1 P_1Q_1 + 2\beta_1 P_2Q_2 + \gamma_1(P_1Q_2 + P_2Q_1) \right] {f} + \hfill   \\

\qquad \qquad \qquad \qquad  + 2\big(\widehat{\Gamma^A}(\mathbf{Q}-\frac{\mathbf{P}}2) - \widehat{\Gamma^A}(\mathbf{Q}+\frac{\mathbf{P}}2) \big)\left[
\xi_1 \int f(\mathbf{P},s,t)ds + \zeta_1
 \int g(\mathbf{P},s,t)ds \right]= 0, \\ { } \\

i \frac{d}{dt} {g} -4\pi^2 \left[ 2\alpha_2 P_1Q_1 + 2\beta_2 P_2Q_2 + \gamma_2(P_1Q_2 + P_2Q_1) \right] {g} + \hfill \\

\hfill + 2\big(\widehat{\Gamma^B}(\mathbf{Q}-\frac{\mathbf{P}}2) - \widehat{\Gamma^B}(\mathbf{Q}+\frac{\mathbf{P}}2) \big)\left[
\xi_2 \int g(\mathbf{P},s,t)ds + \zeta_2
 \int f(\mathbf{P},s,t)ds \right]= 0
\ea
\ee
where of course $\widehat{\Gamma^A}(\mathbf{Q}) = \mathcal{F}_{\mathbf{q}\to\mathbf{Q}}[\Gamma^A(\mathbf{q})]$ and similarly for $\widehat{\Gamma^B}(\mathbf{Q}).$ Observe moreover that
\beq
 \int f(\mathbf{P},s,t)ds = \mathcal{F}^{-1}_{\mathbf{p}\to \mathbf{P}} [n_A(\mathbf{p},t)], \qquad 
  \int g(\mathbf{P},s,t)ds = \mathcal{F}^{-1}_{\mathbf{p}\to \mathbf{P}} [n_B(\mathbf{p},t)],
\ee
hence the notation
\beq
\widecheck{n}_A(\mathbf{P},t) :=  \int f(\mathbf{P},s,t)ds, \qquad 
\widecheck{n}_B(\mathbf{P},t) :=  \int g(\mathbf{P},s,t)ds
\ee
is natural. Also in the interest of brevity let us introduce the bilinear forms
\beq\label{eq:biliform}
\llparenthesis \mathbf{P},\mathbf{Q} \rrparenthesis_j:= 4\pi^2 \left[ 2\alpha_j P_1Q_1 + 2\beta_j P_2Q_2 + \gamma_j(P_1Q_2 + P_2Q_1) \right], \qquad j \in \{1,2\}. 
\ee
Now taking the Laplace transform in time 
\beq
\bac
\widetilde{f}(\mathbf{P},\mathbf{Q},\omega) = \mathcal{L}_{t\to \omega}[{f}(\mathbf{P},\mathbf{Q},t)], \qquad
\widetilde{g}(\mathbf{P},\mathbf{Q},\omega) = \mathcal{L}_{t\to \omega}[{g}(\mathbf{P},\mathbf{Q},t)], \\

\widetilde{n}_A(\mathbf{P},\omega) = \mathcal{L}_{t\to \omega}[\widecheck{n}_A(\mathbf{P},t)], \qquad \widetilde{n}_B(\mathbf{P},\omega) = \mathcal{L}_{t\to \omega}[\widecheck{n}_B(\mathbf{P},t)]
\ea
\ee
we obtain the system
\beq
\bac
i \omega \widetilde{f} -\llparenthesis\mathbf{P},\mathbf{Q}\rrparenthesis_1 \widetilde{f} + 
 2\big(\widehat{\Gamma^A}(\mathbf{Q}-\frac{\mathbf{P}}2) - \widehat{\Gamma^A}(\mathbf{Q}+\frac{\mathbf{P}}2) \big)
\left[
\xi_1  \widetilde{n}_A(\mathbf{P},\omega) + \zeta_1
\widetilde{n}_B(\mathbf{P},\omega) 
\right]= f(\mathbf{P},\mathbf{Q},0), \\ { } \\

i \omega \widetilde{g} -\llparenthesis\mathbf{P},\mathbf{Q}\rrparenthesis_2 \widetilde{g} +%
 2
\big(\widehat{\Gamma^B}(\mathbf{Q}-\frac{\mathbf{P}}2) - \widehat{\Gamma^B}(\mathbf{Q}+\frac{\mathbf{P}}2) \big)
\left[
\xi_2 \widetilde{n}_B(\mathbf{P},\omega) + \zeta_2
 \widetilde{n}_A(\mathbf{P},\omega) 
 \right]= g(\mathbf{P},\mathbf{Q},0).
\ea
\ee
By dividing with the free-space part we obtain
\beq\label{eq:alapltransihgf3243}
\bac
 \widetilde{f}    + 2\frac{\widehat{\Gamma^A}(\mathbf{Q}-\frac{\mathbf{P}}2) - \widehat{\Gamma^A}(\mathbf{Q}+\frac{\mathbf{P}}2) }{i \omega-\llparenthesis\mathbf{P},\mathbf{Q}\rrparenthesis_1 }
 \left[
\xi_1  \widetilde{n}_A + \zeta_1
\widetilde{n}_B 
 \right]=  \frac{f(\mathbf{P},\mathbf{Q},0)}{i \omega-\llparenthesis\mathbf{P},\mathbf{Q}\rrparenthesis_1}, \\ { } \\

\widetilde{g} + 2\frac{\widehat{\Gamma^B}(\mathbf{Q}-\frac{\mathbf{P}}2) - \widehat{\Gamma^B}(\mathbf{Q}+\frac{\mathbf{P}}2) }{
i \omega-4\pi^2 \llparenthesis\mathbf{P},\mathbf{Q}\rrparenthesis_2}
\left[
\xi_2 \widetilde{n}_B + \zeta_2
 \widetilde{n}_A 
 \right]=\frac{g(\mathbf{P},\mathbf{Q},0)}{
i \omega-\llparenthesis\mathbf{P},\mathbf{Q}\rrparenthesis_2}.
\ea
\ee
The two key observations here are the following:
\begin{itemize}
	\item the rhs of both equations correspond to the linear free-space solutions (and thus can be treated as known and well-behaved functions), moreover
	\item we can now integrate both equations in the $\mathbf{Q}$ variables and {obtain a closed system for the position densities} $\widetilde{n}_A,\widetilde{n}_B;$ {achieving this is in fact the motivation for all the transforms and changes of variables}.
\end{itemize}
\beq\label{eq:aformofdthesyst432}
\bac
 \widetilde{n}_A    +  \left[
\xi_1  \widetilde{n}_A + \zeta_1
\widetilde{n}_B 
 \right]  2\int\limits_{\mathbf{Q}}\frac{\widehat{\Gamma^A}(\mathbf{Q}-\frac{\mathbf{P}}2) - \widehat{\Gamma^A}(\mathbf{Q}+\frac{\mathbf{P}}2) }{
i \omega-\llparenthesis\mathbf{P},\mathbf{Q}\rrparenthesis_1}d\mathbf{Q}
=   
\int\limits_{\mathbf{Q}}\frac{f(\mathbf{P},\mathbf{Q},0)}{i \omega-\llparenthesis\mathbf{P},\mathbf{Q}\rrparenthesis_1} d\mathbf{Q}, \\ { } \\

\widetilde{n}_B + \left[
\xi_2 \widetilde{n}_B + \zeta_2
 \widetilde{n}_A 
 \right] 2\int\limits_{\mathbf{Q}}\frac{\widehat{\Gamma^B}(\mathbf{Q}-\frac{\mathbf{P}}2) - \widehat{\Gamma^B}(\mathbf{Q}+\frac{\mathbf{P}}2) }{i \omega-\llparenthesis\mathbf{P},\mathbf{Q}\rrparenthesis_2} d\mathbf{Q}
=\int\limits_{\mathbf{Q}}\frac{g(\mathbf{P},\mathbf{Q},0)}{i \omega-\llparenthesis\mathbf{P},\mathbf{Q}\rrparenthesis_2}d\mathbf{Q}.
\ea
\ee
We can summarize the result by denoting
\beq\label{eq:techdefs}
\bac
\widetilde{n}_{A}^0(\mathbf{P},\omega) := \int\limits_{\mathbf{Q}}\frac{f(\mathbf{P},\mathbf{Q},0)}{i \omega-\llparenthesis\mathbf{P},\mathbf{Q}\rrparenthesis_1} d\mathbf{Q},  \qquad 

\widetilde{n}_{B}^0(\mathbf{P},\omega) := \int\limits_{\mathbf{Q}}\frac{g(\mathbf{P},\mathbf{Q},0)}{i \omega-\llparenthesis\mathbf{P},\mathbf{Q}\rrparenthesis_2}d\mathbf{Q}, \\

h^A(\mathbf{P},\omega):= 2\int\limits_{\mathbf{Q}}\frac{\widehat{\Gamma^A}(\mathbf{Q}+\frac{\mathbf{P}}2) - \widehat{\Gamma^A}(\mathbf{Q}-\frac{\mathbf{P}}2) }{i \omega-\llparenthesis\mathbf{P},\mathbf{Q}\rrparenthesis_1}d\mathbf{Q}, \qquad 

h^B(\mathbf{P},\omega):= 2\int\limits_{\mathbf{Q}}\frac{\widehat{\Gamma^B}(\mathbf{Q}+\frac{\mathbf{P}}2) - \widehat{\Gamma^B}(\mathbf{Q}-\frac{\mathbf{P}}2) }{i \omega-\llparenthesis\mathbf{P},\mathbf{Q}\rrparenthesis_2} d\mathbf{Q}
\ea
\ee
so that  system \eqref{eq:aformofdthesyst432} becomes
\beq
%
\left\{
\bac
\big( 1-\xi_1 h^A \big) \widetilde{n}_{A}  -\zeta_1 h^A \widetilde{n}_{B}= \widetilde{n}_{A}^0, \\ { } \\
 
-\zeta_2 h^B \widetilde{n}_A + (1-\xi_2 h^B) \widetilde{n}_B = \widetilde{n}_{B}^0
\ea
\right\}
\ee
leading to
\beq
\bac
\widetilde{n}_A = 
\frac{1-\xi_2 h^B}{(1-\xi_1 h^A)(1-\xi_2 h^B) - \zeta_1\zeta_2 h^Ah^B} \widetilde{n}_A^0 + \frac{\zeta_1 h^A}{(1-\xi_1 h^A)(1-\xi_2 h^B) - \zeta_1\zeta_2 h^Ah^B}\widetilde{n}_B^0,\\
\widetilde{n}_B = 
\frac{-\zeta_2 h^B}{(1-\xi_1 h^A)(1-\xi_2 h^B) - \zeta_1\zeta_2 h^Ah^B} \widetilde{n}_A^0 - \frac{1-\xi_1 h^A}{(1-\xi_1 h^A)(1-\xi_2 h^B) - \zeta_1\zeta_2 h^Ah^B}\widetilde{n}_B^0
\ea
\ee
as long as the  determinant (itself a function of $\mathbf{P},$ $\omega$)  is nonzero,
\[
(1-\xi_1 h^A(\mathbf{P},\omega))(1-\xi_2 h^B(\mathbf{P},\omega)) - \zeta_1\zeta_2 h^A(\mathbf{P},\omega)h^B(\mathbf{P},\omega) \neq 0.
\]
So the problem is altogether linearly stable if eqs. \eqref{eq:stabilitycond8760}, \eqref{eq:stabilitycond8765} hold.

If both conditions \eqref{eq:stabilitycond8760}, \eqref{eq:stabilitycond8765} hold, then
\beq
\widecheck{n}_A(\mathbf{P},t) = \mathcal{L}^{-1} [\widetilde{n}_A(\mathbf{P},\omega)], \quad
\widecheck{n}_B(\mathbf{P},t) = \mathcal{L}^{-1} [\widetilde{n}_B(\mathbf{P},\omega)]
\ee
inherit an $L^2-$type decay in time, and one expects that Landau damping estimates similar to what was done in the scalar case in \cite{Athanassoulis2018} are possible. Even without the rigorous estimates one readily sees that there is no exponential growth possible, i.e. no modulation instability.

%
%
%
%
%
%
%

%
%
\section{Conclusions}\label{sec:conclusions}



%
Using the same approach as in the classical Alber equation, we derived a system applicable to crossing seas, i.e. two quasi-unidirectional wave systems meeting in the ocean. This is a genuinely 2-dimensional situation, and we produce for the first time the stability condition that controls whether Landau damping or modulation instability is present. This enables for the first time a  systematic detection of  modulation instability in crossing seas, including with realistic data. Such a study is now made possible, but  involves several tedious steps that are beyond the scope of the current paper.

In this paper we also study a collection of non-parametric unidirectional spectra; this can be thought of as a first step towards the 2-dimensional study. The first subtlety that becomes apparent is that the choice of $k_0,$ the carrier wavenumber, is an important question that can affect the results in a non-negligible way. (For example, the mean, mode or median wavenumbers are all plausible choices, and we did see spectra where these were significantly different, and would even lead to different stability/instability classifications for a handful of spectra.) Here we used  the $k_0$ provided by the Norwegian Meteorological Institute,
and found that the vast majority of spectra in our collection are stable; however among the most extreme sea-states it is possible to find spectra that would be classified as (barely) unstable. This is in agreement with the main findings of \cite{Athanassoulis2017,Athanassoulis2018}, and indicates that the Penrose-Alber instability really is a limiting factor of how narrow a spectrum can be in the ocean. The fact that there exist spectra all the way up to instability places additional emphasis on the question of how would unstable (or even borderline stable) spectra behave in the field, something that is really not well understood at all.

Moreover, we find that the novel Proximity to Instability (PTI) metric, defined in equation \eqref{def:pti}, correlates well ($>85\%$ Spearman rank correlation) with the Benjamin-Feir Index (BFI) and the steepness of a sea state, as well as with Monte Carlo estimates for the kurtosis and the probability of rogue waves. This validates the PTI as a meaningful metric for the study of sea states in the unidirectional setting, where things are rather well understood. 

Given that now we have a version of the Alber equation (and the corresponding stability condition) applicable to crossing seas, the natural next step is to investigate crossing seas situations and see how well a generalized version of the PTI  metric would correlate e.g. with Monte Carlo estimates of the probability of rogue waves appearing there, as well as directional generalizations of the BFI \cite{Waseda2009,Mori2011}. Another option for this kind of analysis would be by working with the stability condition for the  broadband 
  Crawford-Saffman-Yuen equation (CSY) \cite{Crawford1980,Andrade2020b}. 

\bigskip 

\noindent {\bf Acknowledgement:}
	We thank the Norwegian Meteorological Institute for providing the wave model spectra used in this study.

\bibliographystyle{siam}
\bibliography{notes_biblio_3.bib}

\appendix

\section{Background}

\begin{theorem}[A complex Isselris theorem] \label{thrm:ComplexIsserlis} Following \cite{Reed1962}; see also \cite{Miller1968}.
Let $z(x)$ be a Gaussian, zero-mean, stationary process with the additional property that
\beq\label{eq:circsymm1}
E[u(x)u(x')]=0 \qquad  \forall x,x'\in\R.
\ee
Then 
\[
E[\overline{z(x_1)z(x_2)}z(x_3)z(x_4)]=E[\overline{z(x_1)}z(x_3)]E[\overline{z(x_2)}z(x_4)]+E[\overline{z(x_2)}z(x_3)]E[\overline{z(x_1)}z(x_4)].
\]
\end{theorem}

\noindent {\bf Remark: }
This result directly implies the closure relation
\beq\label{eq:AlbMomCLo}
E[\overline{u(\alpha,t)u(\beta,t)}u(\alpha,t)u(\alpha,t)]=2E[\overline{u(\alpha,t)}u(\alpha,t)]E[\overline{u(\beta,t)}u(\alpha,t)],
\ee
which is exactly equation \eqref{eq:gmcl098}.
Moreover, the condition \eqref{eq:circsymm1} is equivalent to {circular symmetry}, i.e. to the condition that
\beq\label{eq:gaugeinv}
\{ e^{i\theta} u(x)\}_{\theta\in[0,2\pi)} \mbox{ are identically distributed for all } \theta \in[0,2\pi)
\ee
 which is natural for ocean waves.

\end{document}